\documentclass[12pt]{article}
\usepackage{amsmath,amssymb,amsthm,latexsym}
\usepackage[dvips]{graphics}
\usepackage{epsfig, rotate}
\usepackage[english]{babel}
\usepackage{amsfonts}
\usepackage{array}
\usepackage{color}

\newtheorem{theorem}{{Theorem}}

\newtheorem{lemma}{{Lemma}}

\numberwithin{equation}{section} \numberwithin{theorem}{section}
\numberwithin{lemma}{section}\numberwithin{corollary}{section}

\usepackage
{graphics}
\usepackage
{graphicx}

\newcommand{\be}{\begin{equation}}
\newcommand{\ee}{\end{equation}}
\newcommand{\beaa}{\begin{eqnarray*}}
\newcommand{\eeaa}{\end{eqnarray*}}
\newcommand{\bea}{\begin{eqnarray}}
\newcommand{\eea}{\end{eqnarray}}

\newcommand{\bei}{\begin{itemize}}
\newcommand{\eei}{\end{itemize}}

\newcommand{\bd}{\bold}

\numberwithin{equation}{section}

\def\mymu{\mbox{\boldmath$\mu$}}
\def\mynu{\mbox{\boldmath$\nu$}}
\def\E{\mathbb{E}}

\begin{document}

\noindent {\bf \Large  Empirical Distribution of Scaled Eigenvalues for Product of Matrices from the Spherical Ensemble}

\vspace{20pt}
\noindent{\bf Shuhua Chang$^1$, Yongcheng Qi$^2$}

\vspace{10pt}

{\small
\noindent $^1$Coordinated Innovation Center for Computable Modeling in Management Science, Tianjin University of Finance and Economics, Tianjin 300222, PR China. \\ Email: szhang@tjufe.edu.cn

\vspace{10pt}

\noindent $^2$Department of Mathematics and Statistics, University of Minnesota Duluth,
1117 University Drive, Duluth, MN 55812, USA.
\\Email: yqi@d.umn.edu

\date{\today}

\vspace{20pt}

\noindent{\bf Abstract.}
Consider the product of $m$ independent $n\times n$ random matrices from the spherical ensemble for $m\ge 1$.
The empirical distribution based on the $n$ eigenvalues of the product is called the empirical spectral distribution.
Two recent papers by G\"otze, K\"osters and Tikhomirov (2015) and Zeng (2016) obtain the limit of the empirical spectral distribution for the product when $m$ is a fixed integer.
In this paper, we investigate the limiting empirical distribution of scaled eigenvalues for the product of $m$ independent matrices from the spherical ensemble in the case when $m$ changes with $n$, that is, $m=m_n$ is an arbitrary sequence of positive integers.

\vspace{20pt}

}


\newpage

\section{Introduction}\label{intro}

The study  of random matrices has attracted much attention from mathematics and physics communities and has found
applications in areas such as heavy-nuclei (Wigner, 1955),
number theory (Mezzadri and Snaith, 2005), condensed matter physics (Beenakker, 1997),
wireless communications (Couillet and Debbah, 2011), and high dimensional statistics  (Johnstone (2001, 2008) and Jiang (2009)), just to mention a few. We refer the interested reader to the Oxford Handbook of Random Matrix Theory edited by Akemann, Baik and Francesco (2011) for more references and applications in mathematics and physics.

Some recent research focuses on product of random matrices which have applications
in wireless telecommunication, disordered spin chain, the stability of large complex system,
quantum transport in disordered wires,  symplectic maps and Hamiltonian mechanics,
quantum chromo-dynamics at non-zero chemical potential. See, e.g., Ipsen (2015) for details.

Assume that $m\ge 1$ is an integer. Let $\bd{X}_1, \cdots, \bd{X}_m$ be $m$ independent and identically distributed $n\times n$ random matrices. The product of the $m$ matrices
\begin{equation}\label{Xm}
\bd{X}^{(m)}=\bd{X}_1\bd{X}_2\cdots\bd{X}_m
\end{equation}
is an $n\times n$ random matrix.  The limits of the empirical spectral distributions for the product $\bd{X}^{(m)}$ have been studied in the literature. Several authors, e.g., G\"{o}tze and Tikhomirov (2010),
Bordenave (2011), O'Rourke and Soshnikov (2011) and O'Rourke {\it et al.} (2015) have investigated the limiting empirical spectral distribution
for the product from the complex Ginibre ensemble when $m$ is fixed.  G\"otze,  K\"osters and  Tikhomirov (2015) and Zeng (2016) have obtained the limiting empirical spectral distribution for the product from the spherical ensemble when $m$ is fixed.
Jiang and Qi (2015b) have investigated the limiting empirical distribution for eigenvalues of $\bd{X}^{(m)}$ by allowing that $m$ changes with $n$. Jiang and Qi (2015b) also consider the product of truncations of $m$ independent Haar unitary matrices when $m=m_n$ depends on $n$.

In this paper, we consider the product of $m$ independent matrices from the spherical ensemble.
Let $\bd{A}$ and $\bd{B}$ be two $n\times n$ matrices and all of the $2n^2$ entries of the matrices are i.i.d. standard complex normal random variables. Then, $\bd{X}:=\bd{A}^{-1}\bd{B}$ is called a spherical ensemble (Hough {\it et al.}, 2009).
Let $\bd{z}_1, \cdots, \bd{z}_n$  be the eigenvalues of $\bd{X}$. Then, their joint probability density function is given by
\bea\label{elastic}
C_1\cdot\prod_{j<k}|z_j-z_k|^2\cdot \prod_{k=1}^n\frac{1}{(1+|z_k|^2)^{n+1}},
\eea
where $C_1$ is a normalizing constant; see, for example,  Krishnapur (2009).

Let $\bd{X}_1, \cdots, \bd{X}_m$ be $m$ independent and identically distributed $n\times n$ random matrices from the spherical ensemble, that is, they have the same distribution as $\bd{X}$ defined above.  Define the product ensemble $\bd{X}^{(m)}$ as in \eqref{Xm}.
Again, let $\bd{z}_1, \cdots, \bd{z}_n$ be the eigenvalues of $\bd{X}^{(m)}$. Then their joint probability density function is given by
\bea\label{elastic-m}
C_m\cdot\prod_{j<k}|z_j-z_k|^2\cdot \prod_{k=1}^nw_m(z_k),
\eea
where $C_m$ is a normalizing constant, $w_m(z)$ is given by
\[
w_n(z)=\frac{\pi^{m-1}}{(n!)^m}G^{m,m}_{m,m}\left({{(-n,-n,\cdots, -n)_m}\atop{(0,0\cdots, 0)_m}} \Big| |z|^2\right),
\]
and $G^{m,m}_{m,m}\left({{(-n,-n,\cdots, -n)_m}\atop{(0,0\cdots, 0)_m}} \Big| |z|^2\right)$ is a Meijer $G$-function. See Adhikari {\it et al.} (2016). A recursive formula for $w_m$ is given by
\[
w_{k+1}(z)=2\pi\int^\infty_0w_k(\frac{z}{r})\frac{1}{(1+r^2)^{n+1}}\frac{d\,r}{r}
\]
for $k\ge 1$ with initial $w_1(z)=\displaystyle\frac{1}{(1+|z|^2)^{n+1}}$, which is obtained by Zeng (2016). Obviously, \eqref{elastic-m} reduces to \eqref{elastic} when $m=1$.

Define the empirical spectral distribution (or measure)
\begin{equation}\label{ESD-m}
\mymu_n^*=\frac{1}{n}\sum^n_{j=1}\delta_{\bd{z}_j}.
\end{equation}
If $m\ge 1$ is a fixed integer, Zeng (2016) has proved that
\begin{equation}\label{mun*}
\mymu_n^*\mbox { converges weakly to a distribution } \mymu^* \mbox { with a density function } p_m(z)
\end{equation}
with probability one,  where $p_m(z)$ is given by
\begin{equation}\label{measure-m}
p_m(z)=\frac{1}{m\pi}\frac{|z|^{2/m-2}}{(1+|z|^{2/m})^2}, ~~~z\in \mathbb{C},
\end{equation}
and $\mathbb{C}$ denotes the complex plane.  The universality of convergence in \eqref{mun*} is proved
by G\"otze, K\"osters and  Tikhomirov (2015). More precisely, without assuming the normality,
G\"otze, K\"osters and  Tikhomirov (2015) show that
\eqref{mun*} holds in probability for a large class of random matrices satisfying the Lindeberg condition.

\vspace{10pt}
For the spherical ensemble i.e. $m=1$, \eqref{mun*} has been proved in Bordenave (2011). In fact, Bordenave (2011) has obtained a universal result for the spherical ensemble without assuming the normality of entries in random matrices $\bd{A}$ and $\bd{B}$. The maximum absolute value of the eigenvalues, $\max_{1\le j\le n}|\bd{z}_j|$, is called the spectral radius.
For the spherical ensemble, the limiting distribution for the spectral radius has been obtained in Jiang and Qi (2015a).

In this paper, we will assume that $\{m_n, ~n\ge 1\}$ is an arbitrary sequence of positive
integers and consider the product of $m_n$ independent matrices from the spherical ensemble.  We are interested
in the limiting empirical spectral distribution of the product ensemble $\bd{X}^{(m_n)}$. By defining a new empirical measure based on properly scaled eigenvalues of the product ensemble, we show that the limiting empirical distribution exists
and is free of the sequence $\{m_n\}$.
In particular, our result can reduce to \eqref{mun*} when $m_n=m$, where $m\ge 1$ is any fixed integer.

\section{Main Result}\label{main}

As we assume that $m_n$ can change with $n$, our goal is to  define
the empirical spectral distribution in a different way than \eqref{ESD-m}
so that the limiting distribution is free of the sequence $\{m_n\}$.
Note that the eigenvalues $\bd{z}_1, \cdots, \bd{z}_n$ for the product $\bd{X}^{(m_n)}$ defined in \eqref{Xm} are complex random variables.
Write
\[
\boldsymbol{\theta}_j=\arg(\bd{z}_j)\in [0, 2\pi)\ \ \ \mbox{such that}\ \ \ \bd{z}_j=|\bd{z}_j|\cdot e^{i\boldsymbol{\theta}_j}
\]
for $1\le j\le n$.
To achieve our goal, we define the empirical distribution based on scaled eigenvalues as
\begin{equation}\label{mun}
\mymu_n=\frac{1}{n}\sum^n_{j=1}\delta_{(\boldsymbol{\theta}_j, |\bd{z}_j|^{1/m_n})}.
\end{equation}

We have the following result on the convergence of $\mymu_n$.

\begin{theorem}\label{thm1} With probability one,  $\mymu_n$ converges weakly to a probability measure $\mymu$ with density
\begin{equation}\label{density}
f(\theta, r)=\frac{1}{\pi}\frac{r}{(1+r^2)^2}, ~~~\theta\in [0, 2\pi), r\in (0, \infty).
\end{equation}
\end{theorem}

\vspace{10pt}

\noindent{\bf Remark 1.} A complex number $z=re^{i\theta}$ should be interpreted as a $2$-dimensional vector $(r\cos(\theta), r\sin(\theta))$ in the definition of the empirical spectral distribution given in \eqref{ESD-m}.  Now consider transformation
$z=\xi(\theta, r)=re^{i\theta}$, $\theta\in [0, 2\pi), r\in (0, \infty)$. The Jacobian for this transformation is equal
to $r=|z|$. Therefore, if we assume that $(\bd{r}, \boldsymbol{\theta})$ is a random vector with probability density $f(\theta,r)$ given in \eqref{density}, then the density function for $\bd{z}=\bd{r}e^{i\boldsymbol{\theta}}$ is
\begin{equation}\label{ff}
f(z)=\frac{1}{\pi}\frac{|z|}{(1+|z|^2)^2}\frac{1}{|z|}=\frac{1}{\pi}\frac{1}{(1+|z|^2)^2}, ~~z\in \mathbb{C}.
\end{equation}
Now we can apply the continuous mapping theorem and restate Theorem~\ref{thm1}
as follows: with probability one, the empirical distribution
\[
\frac{1}{n}\sum^n_{j=1}\delta_{|\bd{z}_j|^{1/m_n}e^{i\boldsymbol{\theta}_j}}
=\frac{1}{n}\sum^n_{j=1}\delta_{\xi(\boldsymbol{\theta}_j, |\bd{z}_j|^{1/m_n})}=\mymu_n\circ\xi^{-1}
\]
converges weakly to a probability distribution $\mymu\circ\xi^{-1}$ which has density function $f(z)$ defined in \eqref{ff}.

\noindent{\bf Remark 2.} When $m_n=m$ for all $n$, where $m\ge 1$ is a fixed integer, we can show that Theorem~\ref{thm1} implies \eqref{mun*}. In fact, this can be seen from a simple transformation:   $z=\xi(\theta, r)=r^{m}e^{i\theta}$. The Jacobian for this transformation is $mr^{2m-1}=m|z|^{2-1/m}$.
Again, as in Remark 1, if we assume that $(\bd{r}, \boldsymbol{\theta})$ is a random vector with the probability density given in \eqref{density}, then the density function for $\bd{z}=\bd{r}^{m}e^{i\boldsymbol{\theta}}$ is
\[
\frac{1}{\pi}\frac{|z|^{1/m}}{(1+|z|^{2/m})^2}\frac{1}{m|z|^{2-1/m}}=\frac{1}{m\pi}\frac{|z|^{2/m-2}}{(1+|z|^{2/m})^2}, ~~z\in \mathbb{C},
\]
which is the same as $p_m(z)$ defined in \eqref{measure-m}. Now we can apply the continuous mapping theorem and obtain that with probability one, the empirical distribution
\[
\mymu_n^*=\frac{1}{n}\sum^n_{j=1}\delta_{\bd{z}_j}
=\frac{1}{n}\sum^n_{j=1}\delta_{\xi(\boldsymbol{\theta}_j, |\bd{z}_j|^{1/m})}=\mymu_n\circ\xi^{-1}
\]
converges weakly to a probability distribution $\mymu\circ\xi^{-1}$ which has density function $p_m(z)$.

\noindent{\bf Remark 3.} Eigenvalues with a joint density with a similar structure to \eqref{elastic-m} form a
determinantal point process. See, e.g.,  Hough {\it et al.} (2009) for properties of determinantal point processes.
Eigenvalues from the product of Ginibre ensembles and the product of truncations of independent
Haar unitary matrices can be also modeled by determinantal point processes. By developing a special technique for determinantal point processes, Jiang and Qi (2015b) have obtained the limits for the empirical spectral distributions for the two aforementioned product ensembles.

\section{Proof}\label{proof}

The proof of the theorem relies on applications of Theorem 1 and Lemma 2.1 in Jiang and Qi (2015b).

Let $Y_1, \cdots, Y_n$ be $n$ independent positive random variables such that the density function of $Y_j$ is proportional to
$y^{2j-1}w_{m_n}(y)I(y>0)$ for $1\le j\le n$, where $I(A)$ denotes the indictor function of a measurable set $A$.  Define the empirical distribution of $Y_1, \cdots, Y_n$  as
\[
\mynu_n=\frac{1}{n}\sum^n_{j=1}\delta_{Y_j^{1/m_n}}.
\]



Assume that $\{s_{j, \ell}, ~1\le \ell\le m_n,~1\le j\le n\}$ are independent random variables, and the density of $s_{j,\ell}$ is proportional to
$\frac{y^{j-1}}{(1+y)^{n+1}}I(y>0)$ for any $1\le \ell\le m_n,  1\le j\le n$.

Let $\mbox{Unif}[0,2\pi]$ denote the uniform distribution over $[0,2\pi]$ and $\mynu$ denote the probability measure defined on $(0,\infty)$ with
density function $\frac{2r}{(1+r^2)^2}$, $r>0$. Then we see that the probability measure $\mymu$ with density $f(\theta, r)$ given in \eqref{density} is the product measure  of two probability measures $\mbox{Unif}[0,2\pi]$ and $\mynu$, that is, $\mymu=\mbox{Unif}[0,2\pi]\otimes\mynu$.

We have the following conclusions
in our special situation in the present paper.

\noindent{\bf Result 1.} If $\mynu_n$ converges weakly to $\mynu$ with probability one, then
$\mymu_n$ converges weakly to $\mymu$ with probability one. See Theorem 1 in Jiang and Qi (2015b).

\noindent{\bf Result 2.} If for every $r>0$
\begin{equation}\label{Ylimit}
G_n(r):=\frac1n\sum^n_{j=1}P(Y_j^{1/m_n}\le r)\to \frac{r^2}{1+r^2} ~~\mbox{ as } n\to\infty,
\end{equation}
then $\mynu_n$ converges weakly to $\mynu$ with probability one. Note that the limit $\frac{r^2}{1+r^2}$ in \eqref{Ylimit} is equal to  $\mynu((0; r])$.   See Lemma 2.1 in Jiang and Qi (2015b).

Therefore, to complete the proof of the theorem, it suffices to show \eqref{Ylimit}.  To this end, we list some important results we will use in the proof.

\noindent{\bf Result 3.} For each $1\le j\le n$,
\begin{equation}\label{F1}
Y_j^2 \mbox{ and } \prod^{m_n}_{\ell=1}s_{j,\ell} \mbox{ are identically distributed.}
\end{equation}
See Lemma 2.1 in Zeng (2016).

\noindent{\bf Result 4.} For $1\le \ell\le m_n,  1\le j\le n-2$,
\begin{equation}\label{F2:1}
\mu_{j,\ell}:=\E(s_{j,\ell})=\frac{j}{n-j}, ~~~\text{Var}(s_{j,\ell})=\frac{nj}{(n-j)^2(n-j-1)},
\end{equation}
and
\begin{equation}\label{F2:2}
\E\left(\eta(\frac{s_{[nx],\ell}}{\mu_{[nx], \ell}})\right)\to 0~~~\mbox{ as }n\to\infty
\end{equation}
 for any $x\in (0,1)$, where $\eta(y):=y-1-\log (y)\ge 0$ for $y>0$, and $[nx]$ denotes the integer part of $nx$. See the proof of Lemma 2.3 in Zeng (2016).

\noindent{\bf Result 5.} $Y_1^2,\cdots, Y_n^2$ are stochastically increasing, that is,
\begin{equation}\label{F3}
P(Y_1^2\le x)\ge P(Y_2^2\le x)\ge \cdots\ge P(Y_n^2\le x)~~~x\ge 0.
\end{equation}
See Lemma 2.3 in Zeng (2016).

\begin{lemma}\label{lem1} We have
\begin{equation}\label{bridge2}
\log(Y_{[nx]}^{2/m_n})\to \log\frac{x}{1-x}~~\mbox{ in probability}
\end{equation}
for each $x\in (0,1)$.
\end{lemma}

\noindent{\it Proof.} Define $Y_{nj}^2=\prod^{m_n}_{\ell=1}s_{j,\ell}$. We will first show that
\begin{equation}\label{bridge}
\log(Y_{n,[nx]}^{2/m_n})=\frac{1}{m_n}\sum^{m_n}_{\ell=1}\log(s_{[nx],\ell})\to \log\frac{x}{1-x}~~\mbox{ in probability}
\end{equation}
for any $x\in (0,1)$, which is equivalent to
 \begin{equation}\label{bridge1}
\frac{1}{m_n}\sum^{m_n}_{\ell=1}\log\left(\frac{s_{[nx],\ell}}{\mu_{[nx],1}}\right)\to 0~~\mbox{ in probability}
\end{equation}
since $\mu_{[nx],1}=\frac{[nx]}{n-[nx]}\to \frac{x}{1-x}$ from \eqref{F2:1}.

From the definition of $\eta$ given in \textbf{Result 4} we have
\begin{equation}\label{link}
\frac{1}{m_n}\sum^{m_n}_{\ell=1}\log\left(\frac{s_{[nx],\ell}}{\mu_{[nx],1}}\right)
=\frac{1}{m_n}\sum^{m_n}_{\ell=1}\left(\frac{s_{[nx],\ell}}{\mu_{[nx],1}}-1\right)-\frac{1}{m_n}\sum^{m_n}_{\ell=1}
\eta\left(\frac{s_{[nx],\ell}}{\mu_{[nx],1}}\right).
\end{equation}
Since $s_{j, 1}, \cdots, s_{j, m_n}$ are i.i.d. random variables for each $1\le j\le n$, we have from \eqref{F2:1} that
\[
\E\left(\frac{1}{m_n}\sum^{m_n}_{\ell=1}\left(\frac{s_{[nx],\ell}}{\mu_{[nx],1}}-1\right)\right)^2=\frac{1}{m_n}
\frac{\text{Var}(s_{[nx],1})}{\mu_{[nx],1}^2}=\frac{1}{m_n}\frac{n}{[nx](n-[nx]-1)}\to 0
\]
and
\[
\E\left(\frac{1}{m_n}\sum^{m_n}_{\ell=1}\eta\left(\frac{s_{[nx],\ell}}{\mu_{[nx],1}}\right)\right)=\E\left(\eta\left(\frac{s_{[nx],1}}{\mu_{[nx],1}}
\right)\right)\to 0
\]
as $n\to\infty$. From Chebyshev's inequality,
$\frac{1}{m_n}\sum^{m_n}_{\ell=1}(\frac{s_{[nx],\ell}}{\mu_{[nx],1}}-1)$
 converges to zero in probability as $n\to\infty$, and so does
$\frac{1}{m_n}\sum^{m_n}_{\ell=1}\eta(\frac{s_{[nx],\ell}}{\mu_{[nx],1}})$ since $\frac{1}{m_n}\sum^{m_n}_{\ell=1}\eta(\frac{s_{[nx],\ell}}{\mu_{[nx],1}})\ge 0$. In view of \eqref{link}, \eqref{bridge1} is proved and so is \eqref{bridge}. Consequently,  \eqref{bridge2} follows from \eqref{bridge} and \eqref{F1}.  This completes the proof of Lemma~\ref{lem1}.
\hfill$\blacksquare$

Now we turn to prove \eqref{Ylimit}.

Fix a $r\in (0,\infty)$. By setting  $x=\frac{r^2}{1+r^2}$, we have $x\in (0,1)$ and $r^2=\frac{x}{1-x}$.
Now we choose a small $\delta\in (0,1)$
such that $x+\delta\in (0,1)$ and $x-\delta\in (0,1)$. Then it follows from \eqref{bridge2} that
\[
\log(Y_{[n(x+\delta)]}^{2/m_n})\to \log\frac{x+\delta}{1-(x+\delta)}~~\mbox{ in probability}
\]
and
\[
\log(Y_{[n(x-\delta)]}^{2/m_n})\to \log\frac{x-\delta}{1-(x-\delta)}~~\mbox{ in probability}.
\]
 Since $\delta_1:=\log\frac{x+\delta}{1-(x+\delta)}-\log(r^2)>0$ and $\delta_2:=\log(r^2)-\log\frac{x-\delta}{1-(x-\delta)}>0$, we obtain
 from the above two equations that
 \begin{equation}\label{plus}
 P(Y_{[n(x+\delta)]}^{2/m_n}\le r^2)=P\left(\log(Y_{[n(x+\delta)]}^{2/m_n})-\log\frac{x+\delta}{1-(x+\delta)} \le-\delta_1\right)\to 0
 \end{equation}
 and
  \begin{equation}\label{minus}
 P(Y_{[n(x-\delta)]}^{2/m_n}> r^2)=P\left(\log(Y_{[n(x-\delta)]}^{2/m_n})-\log\frac{x-\delta}{1-(x-\delta)} > \delta_2\right)\to 0.
 \end{equation}
Therefore, by using \eqref{F3} and \eqref{plus} we have
\begin{eqnarray*}
G_n(r)
      &=&\frac1n\sum^{[n(x+\delta)]-1}_{j=1}P(Y_j^{2/m_n}\le r^2)+\frac1n\sum^n_{j=[n(x+\delta)]}P(Y_j^{2/m_n}\le r^2)\\
      &\le&\frac{[n(x+\delta)]-1}n+\frac{n-[n(x+\delta)]+1}nP(Y_{[n(x+\delta)]}^{2/m_n}\le r^2)\\
      &\to& x+\delta.
\end{eqnarray*}
Similarly, in view of  \eqref{F3} and \eqref{minus} we obtain
\begin{eqnarray*}
G_n(r)
&=&\frac1n\sum^{[n(x-\delta)]}_{j=1}P(Y_j^{2/m_n}\le r^2)+\frac1n\sum^n_{j=[n(x-\delta)]+1}P(Y_j^{2/m_n}\le r^2)\\
      &\ge&\frac{[n(x-\delta)]-1}nP(Y_{[n(x-\delta)]}^{2/m_n}\le r^2)\\
      &\to& x-\delta.
\end{eqnarray*}
Consequently, we prove that
\[
x-\delta\le \liminf_{n\to\infty}G_n(r)\le \limsup_{n\to\infty}G_n(r)\le x+\delta.
\]
By letting $\delta\to 0$ on both sides above we get $\lim_{n\to\infty}G_n(r)=x=\frac{r^2}{1+r^2}$, that is, \eqref{Ylimit} holds.

This completes the proof of Theorem~\ref{thm1}.

\vspace{20pt}

\noindent{\bf Acknowledgements:}  We would like to thank two reviewers for their constructive suggestions that have led to improvement in the layout and readability of the paper. Chang's research was supported in part by the Major Research Plan of the National Natural Science Foundation of China (91430108), the National Basic Research Program (2012CB955804), the National Natural Science Foundation of China (11171251), and the Major Program of Tianjin University of Finance and Economics (ZD1302).

\baselineskip 12pt
\def\ref{\par\noindent\hangindent 25pt}

\end{document}